\documentclass[jath]{apjrnl}

\equationsection
{\theoremstyle{plain}%
  \newtheorem{theorem}{Theorem}
  \newtheorem{corollary}{Corollary}
  \newtheorem{proposition}{Proposition}
  \newtheorem{lemma}{Lemma}%
}

{\theoremstyle{remark}

\newtheorem{remark}{Remark}
}

{\theoremstyle{definition}
\newtheorem{definition}{Definition}

}

\usepackage{amsfonts}
\usepackage{amsmath}
\usepackage{amssymb}
\usepackage{latexsym}
\usepackage{fancybox}
\usepackage{mathrsfs}
\usepackage{epsfig}

\def\eqref{\eqnref}
\def\b{\mathbb}
\DeclareMathOperator {\si}{si}

\begin{document}

\title{Asymptotic analysis for the Dunkl kernel}

\author{Margit R\"osler
\\Mathematisches Institut, Universit\"at G\"ottingen
\\Bunsenstr. 3--5, D-37073 G\"ottingen, Germany
\\E-mail: roesler@uni-math.gwdg.de
\\and
\\Marcel de Jeu
\\Korteweg-de Vries Institute for Mathematics
\\University of Amsterdam
\\Plantage Muidergracht 24, 1018 TV Amsterdam, The Netherlands
\\E-mail: mdejeu@science.uva.nl
}

\date{February 10, 2002}

\maketitle

\begin{abstract}
This paper studies the asymptotic behavior of the integral kernel of the
Dunkl transform, the so-called Dunkl kernel, when one of its
arguments is fixed and the other tends to infinity either
 within a Weyl chamber
of the associated reflection group, or within a suitable complex domain.
The obtained results are based on the asymptotic analysis of
an associated system of ordinary differential equations.
They generalize the well-known asymptotics of the confluent hypergeometric
function $\phantom{}_1F_1$  to
the higher-dimensional setting and include a complete
short-time asymptotics for the Dunkl-type heat kernel.
As an application, it is shown that
the representing measures of Dunkl's intertwining operator are generically
continuous.
\end{abstract}

\begin{subject}[2000 AMS Subject Classification]
Primary: 33C52; Secondary: 33C67
\end{subject}

\section{Introduction and results}

In the theory of rational Dunkl operators as initiated by
C. F. Dunkl in \cite{D1}, there is an analogue
of the classical exponential function, commonly called the Dunkl kernel.
It generalizes the
 usual exponential function in many respects, and can  be
characterized as the solution of a joint eigenvalue problem for the associated Dunkl operators.
Generally speaking, Dunkl operators are parametrized differential-reflection operators attached
to a finite reflection group. During the last decade, such operators have found considerable
attention in various areas of mathematics and mathematical physics. They are, for example,
useful in the study of integrable quantum many body systems of Calogero-Moser-Sutherland type
(for an up-to-date bibliography, we refer to \cite{DV}), and have led to a rapid development in
the theory of special functions related with root systems; see for instance \cite{DX} and
\cite{H}. Among the variants of Dunkl operators, especially the rational ones allow for
 a far-reaching
harmonic analysis in
close analogy to the classical Fourier analysis on $\b R^N$. For example,
there exists an analogue of the  Fourier transform -
the Dunkl transform -
which  establishes a
natural correspondence between the action of Dunkl operators on one hand and
multiplication operators on the other (\cite{D3}, \cite{dJ1}).
The integral kernel of this transform is the Dunkl kernel.
It was first introduced in \cite{D2} and has since then been
studied in a variety of aspects among which we mention
\cite{D3}, \cite{dJ1}, \cite{O} and \cite{R}.
The present paper
contributes to a  further study of its asymptotic and structural properties.

In order to introduce our setting and results, let $R$ be a reduced (not necessarily
crystallographic) root system in $\b R^N$ with the standard Euclidean inner product
$\langle\,.\,,\,.\,\rangle$. This means that $R\subset\b R^N\setminus\{0\}$ is finite with
$\sigma_\alpha R= R$ and $R\cap \b R\alpha = \{\pm\alpha\}$ for all
 $\alpha\in R,$ with $\sigma_\alpha$ denoting the reflection
in the hyperplane $H_\alpha$ orthogonal to $\alpha$.
Let $G \subset O(\b R^N)$  denote the finite reflection group
generated by the $\sigma_\alpha,\,\alpha\in R$, and put
$\b R^N_{\,reg}:= \b R^N\setminus \cup_{\alpha\in R} H_\alpha$.
The connected components of $\b R^N_{\,reg}$ are called the Weyl
chambers of $G$.
As customary in this context, we  assume that
the root system $R$
is normalized by $|\alpha|^2 = 2$ for all $\alpha,$ and we denote
the bilinear extension of $\langle\,.\,,\,.\,\rangle$ to $\b C^N$
 by $\langle\,.\,,\,.\,\rangle$ as well.  Let further
$k: R\to \b C$ be a multiplicity function
on  $R$ (i.e.\ $k$ is  invariant under the
natural action of $G$ on $R$.) In the present paper we shall assume
throughout that
$k$ is non-negative, i.e. $k(\alpha)\geq 0$ for all $\alpha\in R$.
The (rational) Dunkl operators associated with $G$ and $k$ are given by
\[ T_\xi(k) f(x)  := \partial_\xi f(x) + \sum_{\alpha\in R_+}
   k(\alpha)\langle\alpha,\xi\rangle
\frac{f(x)- f(\sigma_\alpha x)}{\langle\alpha,x\rangle}\,, \quad \xi \in \b R^N.\]
Here $\partial_\xi$ denotes the usual partial derivative in direction
 $\xi$
and $R_+$ is an arbitrary but fixed
positive subsystem of $R$.
It is a remarkable property of the $T_\xi(k)$ that
they mutually commute, see \cite{D1}. The Dunkl kernel $E_k(x,y)$ associated with $G$ and $k$
can be characterized as the unique solution of the joint eigenvalue problem for the corresponding
Dunkl operators, more precisely: for each fixed
 $y\in \b C^N$, the function
$x\mapsto E_k(x,y)$ is  the unique real-analytic solution of the system
\begin{equation}\label{(2.66)}
T_\xi(k) f = \langle\xi,y\rangle f \quad\text{ for all }\xi\in \b R^N
\quad\text{and }\,\quad f(0)=1,
\end{equation}
c.f. \cite{O}. In case $k=0$ we just have $E_k(x,y) = e^{\langle x,y\rangle}$.
The generalized exponential
kernel $E_k(x,y)$ is symmetric
in its arguments and has a unique holomorphic extension to $\b
C^N\times \b C^N$. It satisfies
\begin{equation}\label{(2.71)}
E_k(\lambda z,w) = E_k(z,\lambda w)\, \text{ and } E_k(gz,gw) =
E_k(z,w)\end{equation}
for all $z,w\in \b C^N, \, \lambda\in \b C$ and $g\in G.\,$

Originally, $E_k$ was defined in \cite{D2} by means of the so-called intertwining
operator $V_k$.
In fact, there exists
a unique linear
isomorphism $V_k$ on polynomials, homogeneous of degree $0,$ and such that
\[T_\xi(k)V_k\,=\, V_k\partial_\xi \quad\text{for all }\,
\xi \in \b R^N \,\,\text{ and }\,
V_k(1) = 1; \]
see \cite{D1, DJO}.
 In \cite{R} it was shown that
$V_k$ ($k$ always being non-negative)
admits a positive integral representation as follows: Let $M^1(\b R^N)$ denote the
space of probability measures on the Borel $\sigma$-algebra of $\b R^N$. Then
for every $x\in \b R^N$, there exists a
unique
$\mu_x^k\,\in M^1(\b R^N)$   such that
\begin{equation}\label{(2.2)}
 V_k p(x)\,=\, \int_{\b R^N}p(\xi)\,d\mu_x^k(\xi)
\end{equation}
for each polynomial function $p$ on $\b R^N$.
The representing measures $\mu_x^k$ are compactly
supported with
$\,{\rm supp}\,\mu_x^k \subseteq \, \text{co}\{gx, g\in G\}$, the convex hull of the orbit of
 $x$ under $G$.
By means of formula \eqref{(2.2)}, $V_k$ may be extended  to various larger function spaces,
e.g.  the space $C(\b R^N)$ of continuous functions on
$\b R^N$. We denote this extension  by $V_k$ again. Then for fixed $y\in \b C^N$,
\begin{equation}\label{(2.77)}
 E_k(x,y)= V_k(e^{\langle\,.\,, y\rangle})(x)\,=\, \int_{\b R^N}
e^{\langle \xi, y\rangle}d\mu_x^k(\xi) \quad (x\in \b R^N).
\end{equation}
This in particular implies that
\begin{equation}\label{(2.72)}
E_k(-ix,y) = \overline{E_k(ix,y)}\,\text{ and }\,\, |E_k(ix,y)|\leq 1 \text{ for all }
x,y\in \b R^N.
\end{equation}

As already indicated, the Dunkl kernel is especially of interest as it gives rise
to a corresponding integral transform on $\b R^N$. The Dunkl transform associated
with $G$ and $k$
involves the weight function
\[w_k(x)\,=\, \prod_{\alpha\in R_+} |\langle
\alpha,x\rangle|^{2k(\alpha)},\]
which is $G$-invariant and homogeneous of degree $2\gamma$,
with the index
\[ \gamma:= \gamma(k) = \sum_{\alpha\in R_+} k(\alpha)\,\geq 0.\]
It is defined by
\[\widehat f^{\,k}(\xi):= \, c_k^{-1}\int_{\b R^N} f(x) E_k(-i\xi,x)
\,w_k(x)dx, \quad f\in L^1(\b R^N, w_k);\]
here $c_k$ is the Mehta-type constant
\[ c_k:= \int_{\b R^N} e^{-|x|^2/2} w_k(x)dx.\]
An explicit expression for $c_k$ can be found in \cite{O}. The Dunkl transform shares many
properties of the classical Fourier transform. For example, there exist a Plancherel theorem and
an inversion theorem for it. For details the reader is referred to \cite{D3} and \cite{dJ1}.

In this paper, we study the asymptotic behavior of $x\mapsto E_k(x,y)$
for large  arguments $x$, with
$y\in \b R^N_{reg}$ considered as a fixed parameter.
Let $C$
denote the Weyl chamber attached with the positive subsystem $R_+$,
\[C= \{x\in \b R^N: \langle
\alpha,x\rangle >0 \,\,\text{ for all }\, \alpha\in R_+\},\]
and for $\delta >0,$
\[C_\delta :=
\{x\in C: \langle\alpha,x\rangle\,>\,
\delta |x| \,
  \text{ for all }\,\alpha\in R_+\}.\]
Our main result is the following  asymptotic
behaviour, uniform for the variable tending to infinity in cones $C_\delta$:

\begin{theorem}\label{T:main1}
There exists a constant non-zero vector $v=(v_g)_{g\in G}\in\b C^{|G|}$ such
that for all $y\in C, \,g\in G$ and each $\delta>0,$
\[\lim_{|x|\to\infty,\,x\in C_\delta} {\sqrt{w_k(x)w_k(y)}}\,e^{-i\,\langle x,gy\rangle} E_k(ix,gy)\,=\, v_g.\]
\end{theorem}

\noindent
The proof of this Theorem, which will be given in Section 4, is
based on the analysis of an associated  system of
first order ordinary differential equations, which is derived
from the eigenfunction characterization
 \eqref{(2.66)} of $E_k$.  The idea for this approach goes
 back to \cite{dJ1}, where it was used to obtain
exponential estimates for the Dunkl kernel. An immediate consequence of Theorem 1 is the
following ray asymptotic for the Dunkl kernel, making precise  a conjectural remark in \cite{D3}.

\begin{corollary}\label{C:ray}
For all $x,y\in C$ and $g\in G$,
\[ \lim_{t\to\infty} t^\gamma \,e^{-it\,\langle x,gy\rangle} E_k(itx,gy)\,=\,
   \frac{v_g}{\sqrt{w_k(x)w_k(y)}}\,,\]
the convergence being locally uniform with respect to the parameter $x$.
\end{corollary}

In the particular case $g=e$ (the unit of $G$), this result
can be extended to a larger range of
complex arguments by use of the Phragm\' en-Lindel\"of principle.
We consider  the closed right half plane $H =\{z\in \b C: \text{Re} z
\geq 0\}$, and denote by  $z\mapsto z^\gamma$  the holomorphic branch in
$\b C\setminus
\{x\in \b R: x\leq 0\}$ with $1^\gamma = 1$. We shall prove

\begin{theorem}\label{T:main3} Let $x,y\in C$. Then
\[ \lim_{z\to\infty, z\in H}
z^\gamma e^{-z\langle x,y\rangle} E_k(zx,y)\,=\,\frac{i^\gamma v_e}{\sqrt{w_k(x)w_k(y)}}.
\]
\end{theorem}

An interesting consequence of this result concerns the short-time behavior of the
Dunkl-type heat kernel
\[ \Gamma_k(t,x,y) = \frac{1}{(2t)^{\gamma +N/2}c_k}e^{-(|x|^2 + |y|^2)/4t}
E_k\bigl(\frac{x}{\sqrt{2t}},
\frac{y}{\sqrt{2t}}\bigr)\]
$(x,y\in \b R^N, \,t>0),$
which was first introduced in \cite{R3}. After suitable normalization, the kernel $\Gamma_k$
behaves for short times
like the free Gaussian heat kernel $\,\Gamma_0(t,x,y)\,=\,(4\pi t)^{-N/2}\,
 e^{-|x-y|^2/4t}$, as conjectured in \cite{R2}. More precisely,

\begin{corollary}\label{C:heat} For all $x,y\in C$,
\[\lim_{t\downarrow 0} \frac{\sqrt{w_k(x)\,w_k(y)}\,\,
\Gamma_k(t,x,y)}{\Gamma_0(t,x,y)}\,=\,1\,.\]
\end{corollary}

Indeed, it is immediate from Theorem 2
that the limit above exists for all $x, y\in C$ and is equal to
$i^\gamma v_e c_0/c_k.$  On the other
hand, Theorem 3.3 of \cite{R2} shows -- based on completely different methods --
that the limit exists and equals $1$ for a
 restricted range of arguments $x,y\in C$.
This combination proves Corollary \ref{C:heat} and at the same time implies the value of $v_e$:
\begin{equation}\label{(1.40)}
 v_e \,=\, i^{-\gamma}\frac{c_k}{c_0}.
\end{equation}
\begin{remark}
The explicit determination of the constants $v_g$ with $g\not=e$ (or of $v_e$ without falling
back to \cite{R2}) is an open problem. The proof of 
Theorem \ref{T:main1} yields the
invariance property
\[ v_g = v_{g^{-1}} \quad\text{ for all } g\in G,\]
but additional techniques seem to be necessary to obtain further information.
\end{remark}
The asymptotic result of Theorem \ref{T:main1} also allows
 to deduce at least a certain amount of information about
 the structure of the intertwining operator $V_k$ and its representing measures
 $\mu_x^k$ according to formula \eqref{(2.2)}. These  measures
are explicitly known in very special cases only, namely essentially for the rank-one-case as
well as the symmetric group $S_3$ (\cite{D4}), and very little is known about their general
structure either. In particular, as to the authors' knowledge, no results towards the continuity
properties of the $\mu_x^k$ have been obtained so far. We shall employ a well-known
characterization of continuous measures due to Wiener by means of their (classical)
Fourier-Stieltjes transform, which coincides with the kernel $E_k(-ix,\,.\,)$ in case of
$\mu_x^k$. (Recall that a measure $\mu\in M^1(\b R^N)$ is called continuous, if $\mu(\{x\}) = 0$
for all $x\in\b R^N$.) Theorem \ref{T:main1} gives just sufficient information on the growth of
the Dunkl kernel in order to apply Wiener's criterion. This yields

\begin{theorem}\label{T:main2} If $\gamma >0$, i.e. apart from the
classical Fourier case,
the measure $\mu_x^k$ is continuous
for all $x\in \b R^N_{reg}$.
\end{theorem}

We conjecture that
the measures $\mu_{x}^k$ are even absolutely continuous with respect to Lebesgue measure
for all regular $x$, provided $k$ is such that $\{\alpha\in R\mid k(\alpha)> 0\}$
spans $\b R^N$.
This is in fact true
in the rank-one case, which provides a simple but illustrative example for
our results. A short discussion of this example is given in Section 2. Section 3
contains the asymptotic analysis of the differential equation associated with
the kernel $E_k$, as well as the implications concerning its asymptotic behavior.
These results are the basis for the proofs of  Theorems 1 - 3,
 which are completed in Section 4.

\begin{remark}
We mention that the group invariant counterpart of $E_k$, called
``generalized Bessel function'' in \cite{O},
can be considered as a natural generalization of the usual one-variable Bessel
function, to which it reduces in the rank-one case (see below). For Weyl groups
and certain half-integer multiplicity parameters $k$, generalized Bessel
functions have an interpretation as spherical functions of a Cartan motion group. For
the details concerning this identification we refer
to \cite{dJ2,O}. For such generalized Bessel functions corresponding to the group case, and
with both the geometric and spectral variable in $C$, asymptotic results are derived in
\cite{BC} which are more precise than can be obtained by averaging the results in Theorem
\ref{T:main1}. Estimates on the generalized Bessel function in the group case with spectral
variable in $C$, but with \emph{arbitrary} geometric variable, can be found in \cite{C}. The
methods in \cite{BC,C} however use the presence of additional ambient structure for these
special values of the multiplicity parameters, and therefore do not apply in our case of
general non-negative $k$.
\end{remark}

\section{Example: The rank-one case}

 Let $N=1$. Then
the only choice of $R$
(being reduced) is $R=\{\pm \sqrt 2\}$. Accordingly,
 $G=\{id,\sigma\} \cong  \b Z_2$ with $\sigma(x) = -x$.
The Dunkl operator $T(k) = T_1(k)$ associated with the multiplicity parameter $k\geq 0$ is given
by
\[ T(k)f(x) \,=\, f^\prime(x) + k\,\frac{f(x) - f(-x)}{x}.\]
The corresponding intertwining operator $V_k$ and the kernel
$E_k$ were determined in \cite{D2},
\cite{D3}.
In particular,
\[ E_k^{\b Z_2}(z,w)\,=\, j_{k-1/2}(izw)\,+\,\frac{zw}{2k+1}\,
j_{k+1/2}(izw). \] Here $j_\alpha$ denotes the normalized spherical Bessel function
\[ j_\alpha(z)\,=\,
2^\alpha\Gamma(\alpha+1)\cdot\frac{J_\alpha(z)}{z^\alpha}\,=\,
\Gamma(\alpha+1)\cdot\sum_{n=0}^\infty
\frac{(-1)^n(z/2)^{2n}}{n!\,\Gamma(n+\alpha+1)}\,.\]
The integral representation \eqref{(2.77)} of $E_k^{\b Z_2}$ is
  given by
\begin{align}\label{(2.30)} E_k^{\b Z_2} (z,w)\,=\,&
\frac{\Gamma(k+1/2)}{\Gamma(1/2)\,\Gamma(k)}\int_{-1}^1
e^{tzw}(1-t)^{k-1}(1+t)^k\,dt\,\notag\\
 =\,& e^{zw}\cdot
\phantom{}_1F_1(k,2k+1,-2zw).
\end{align}
Thus for $x\not=0$, the associated representing measure is
\[d\mu_x^k(u)\,=\, \frac{\Gamma(k+1/2)}{\Gamma(1/2)\,\Gamma(k)}\cdot 1_{[-|x|,|x|]}(u)
\frac{1}{|x|} \bigl(1-\frac{u}{x}\bigr)^{k-1}\bigl(1+\frac{u}{x}\bigr)^{k}du,\]
which is absolutely continuous with respect to Lebesgue measure.
Further, recall the well-known asymptotic expansion of Kummer's
  function $\phantom{}_1F_1$ (see e.g. \cite{AS}, (13.5.1)):
\[\lim_{z\to \infty,\,z\in H}z^k\cdot \phantom{}_1F_1(k,2k+1,-2z)\,=\,
  \frac{\Gamma(2k+1)}{2^k \, \Gamma(k+1)}.\]
Thus the constants $v_e$ and $v_\sigma$ in Theorems \ref{T:main1},\ref{T:main3}
are given by
\[
v_e\,=\,\frac{\Gamma(2k+1)}{2^k\,\Gamma(k+1)}\cdot \,i^{-k}
\quad v_\sigma = \frac{\Gamma(2k+1)}{2^k\,\Gamma(k+1)}\cdot \,i^{k} \,.\]
\medskip

\section{Asymptotics of $E_k$ along curves in a Weyl chamber}

For $x,y\in \b R^N$  define
\[ \phi(x,y) = \,\sqrt{w_k(x)w_k(y)}\,
e^{-i\langle x,y\rangle} E_k(ix,y).\]
Observe that $\phi$  is symmetric in its arguments. In this section, we shall study
the asymptotic behavior of $x\mapsto \phi(x,y)$ along curves in a fixed Weyl chamber,
with the second component $y\in \b R^N_{reg}$ being
fixed. In view of the invariance properties of $w_k$
and $E_k$ under the action of $G$, we may restrict ourselves to the case
$x\in C$ (the chamber associated with $R_+$) and $y\in gC$ for some $g\in G$.
Following  an idea of \cite{dJ1}, we introduce an auxiliary
 vector field
  $\,F = (F_g)_{g\in G}$  on $\b R^N\times\b R^N$ by
\[ F_g(x,y)\,:=\, \phi(x,gy).\]
For fixed $y$, we consider $F$ along a differentiable curve $\kappa: (0,\infty)\to C$.
The eigenfunction characterization \eqref{(2.66)}
of $E_k$ then  translates into a first order ordinary differential
equation for $t\mapsto F(\kappa(t),y)$.
Below, we shall determine the asymptotic behavior of its solutions,
provided $\kappa$ is admissible in the following sense:

\begin{definition} A $C^1$-curve\\
\begin{minipage}[t]{.51\linewidth}
$\kappa: (0,\infty)\to C$ is called \emph{admissible}, if
it satisfies the subsequent conditions:
\parskip=-3pt
\begin{enumerate}\itemsep=-2pt
\item[\rm{(1)}] There exists a constant $\delta >0$ such that
   $\kappa(t)\in C_\delta$ for all $t>0.$
\item[\rm{(2)}] $\lim_{t\to\infty} |\kappa(t)| \,=\infty\,$
 and $\kappa^\prime(t)\in C$ for all $t>0$.
\end{enumerate}
\end{minipage}\hfill
\begin{minipage}[t]{.49\linewidth}
\vspace*{-34pt}
\psfull
$$\epsfbox{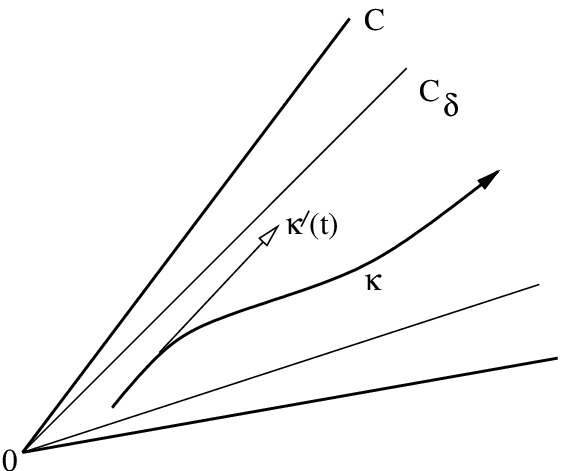}$$
\end{minipage}\hfill
\end{definition}

Notice that conditions (1) and (2) imply that  $\lim_{t\to\infty}\langle\alpha,
\kappa(t)\rangle = \infty$ for all $\alpha\in R_+$.
An important class of admissible curves are the rays $\kappa(t) =tx$ with
some fixed $x\in C$. Corollary \ref{C:ray} describes the asymptotic
behaviour of $x\mapsto F(x,y)$ along such rays.
In the present section, we prove that
 $t\mapsto F(\kappa(t),y)$ is asymptotically constant as
$t\to\infty$ for arbitrary admissible curves, not just for rays. In the next section it will become clear that
the limit value is actually independent of $y$ and $\kappa$.

\begin{theorem} \label{T:Main}  If  $\kappa: (0,\infty)\to C$ is
admissible, then for
every   $y\in C,$  the limit
\[\lim_{t\to\infty} F(\kappa(t),y)\]
exists in $\b C^{|G|}$, and is different from $0$.
\end{theorem}

The subsequent proof of Theorem \ref{T:Main} is based on  the following
variant of formula (3.1) in \cite{dJ1}:

\begin{lemma}\label{L:Kerneq} For fixed $\xi, \,y\in \b R^N,$
\[ \partial_{\xi} F_g(x,y)\,=\, \sum_{\alpha\in R_+} k(\alpha)\,
  \frac{\langle\alpha,\xi\rangle}{\langle\alpha,x\rangle}\,
e^{-i\langle \alpha,\,x\rangle\langle\alpha,\,gy\rangle}\cdot
   F_{\sigma_\alpha g}(x,y)\quad (x\in \b R^N_{reg}).\]
\end{lemma}

\begin{proof} The eigenfunction characterization \eqref{(2.66)}
 of the kernel $E_k$, together
with the invariance property $E_k(gx,gy) = E_k(x,y)$, implies that
\[\partial_{\xi} E_k(x,y)\,=\, \langle\xi,y\rangle E_k(x,y)\,-\,
\sum_{\alpha\in R_+} k(\alpha)
  \frac{\langle\alpha,\xi\rangle}{\langle\alpha,x\rangle}
\bigl( E_k(x,y) - E_k(x,\sigma_\alpha y)\bigr).\]
Moreover, if $x\in \b R^N_{\,reg}$, then
\[ \partial_{\xi} \sqrt{w_k(x)}\,=\, \Bigl(\sum_{\alpha\in R_+} k(\alpha)
\frac{\langle\alpha,\xi\rangle}{\langle\alpha,x\rangle}\Bigr)\cdot
\sqrt{w_k(x)}.\]
It follows that
\begin{align}
 \partial_{\xi}
F_g(x,y)\,&=\,\sqrt{w_k(x)w_k(y)}\, e^{-i\langle x,\,gy\rangle}\cdot
\sum_{\alpha\in R_+} k(\alpha)
   \frac{\langle\alpha,\xi\rangle}{\langle\alpha,x\rangle}\,
   E_k(ix,\sigma_\alpha gy).\notag\\
 &=\,\sum_{\alpha\in R_+} k(\alpha)
   \frac{\langle\alpha,\xi\rangle}{\langle\alpha,x\rangle}\,
 e^{-i\langle x,gy\rangle} e^{i\langle x,\sigma_\alpha gy\rangle}
F_{\sigma_\alpha g}(x,y). \notag
\end{align}
As $\, \sigma_\alpha(x) - x\,=\, -\langle
   \alpha,x\rangle \alpha\,,$
this implies the assertion.
\end{proof}

We are interested in  the derivative of $x\mapsto F(x,y)$ along
differentiable curves in $C$, with the second component $y$ being fixed.
The following is immediate from the previous Lemma:

\begin{corollary}\label{C:Diffsys} For a $C^1$-curve
$\,\kappa: (0,\infty) \to \b R^N_{reg}$
   and fixed $y\in \b R^N$
  define $\, F^\kappa := (F_g^\kappa)_{g\in G}\,$ by $F_g^\kappa(t):=
  F_g(\kappa(t),y).$ Then $F^\kappa$ satisfies the differential equation
\begin{equation}\label{(1.2)}
 (F^\kappa)^\prime(t)\,=\, A^\kappa(t)F^\kappa(t),
\end{equation}
where $\, A^\kappa: (0,\infty) \to \,\b C^{|G|\times |G|}$ is given by
$\,A^\kappa\,=\, \sum_{\alpha\in R_+}
k(\alpha)B_\alpha^\kappa\,$, and
\[(B_\alpha^\kappa(t))_{g,\,h}\,=\,\begin{cases}\displaystyle{
   \frac{\langle\alpha,\kappa^\prime(t)\rangle}{\langle\alpha,\kappa(t)\rangle}
   \,e^{-i\langle\alpha,\kappa(t)\rangle \langle
   \alpha,\,gy\rangle}} & \text{if $h = \sigma_\alpha g$},\\
    0 \,& \text{otherwise.}
    \end{cases}
\]
\end{corollary}

If $\kappa$ is a ray, i.e. $\kappa(t) = tx$ with some fixed $x\in C$, then
\[ (B_\alpha^\kappa(t))_{g,\sigma_\alpha g}\,=\,
   \frac{1}{t}\,e^{-it\langle\alpha,x\rangle\langle\alpha,gy\rangle}.
\]
Notice that in this typical case,  $t=\infty$  is an essential
singularity of $A^\kappa$.
However, $A^\kappa(t) = O\bigl(\frac{1}{t}\bigr)$ as $t\to\infty$, i.e.
the system is asymptotically constant (in fact asymptotically zero) in the sense of \cite{East}.
This suggests that the solution $F^\kappa$ should be asymptotically constant
as $t\to\infty$, which just means that
$\lim_{t\to\infty} F^\kappa(t)$ exists as asserted in Theorem \ref{T:Main}.

 The decisive criterion for the proof of Theorem \ref{T:Main} is the following
result on  the asymptotic integration of ordinary linear differential
equations, which is a special case of the Levinson-type Theorem 1.11.1 in
\cite{East}, and is originally due to \cite{W}:

\begin{proposition}\label{P:Wintner} (\cite{East},\,\cite{W})
Consider the
linear differential equation
\begin{equation}\label{(1.1)}
 x^\prime(t)\,=\, A(t) x(t),
\end{equation}
where  $A:[t_0,\infty)\to \b C^{n\times n}$ is a continuous matrix-valued
function satisfying
the following  integrability conditions:
\parskip=-3pt
\begin{enumerate}\itemsep=-2pt
\item[\rm{(1)}] The matrix-valued improper Riemann
   integral $\,\int_{t_0}^\infty A(t)dt\,$ converges. In particular,
$\,\widetilde A(t):=
\int_t^\infty A(s)ds\,$ is well-defined on $[t_0,\infty)$.
\item[\rm{(2)}] $A\widetilde A \in  L^1\bigl([t_0,\infty), \b C^{\,n\times n}\bigr).$
\end{enumerate}
Then $\eqref{(1.1)}$ has a basis of solutions $\{x_k(t), 1\leq k \leq n\}$ of the
asymptotic form $x_k(t) = e_k + o(1)$ as $t\to \infty$, where $e_k$ is the
$k$-th unit vector in $\b R^n$.  In particular,
for each solution $x$ of $\eqref{(1.1)},$ the limit
   $\, \lim_{t\to\infty} x(t)\,$ exists. Moreover, if
 $x\not= 0$, then $\lim_{t\to\infty}x(t)\not= 0.$
\end{proposition}

\begin{proof}[Proof of Theorem \ref{T:Main}]
We  shall verify that the matrix $A^\kappa$ satisfies the conditions
of Proposition \ref{P:Wintner} with arbitrary $t_0>0$.  For
(1), let $g\in G$ and $\alpha\in R_+$.
Then if $T>t\geq t_0$,
\begin{align} \int_t^T \bigl(B_\alpha^\kappa(s)\bigr)_{g,\sigma_\alpha g}\,  ds \,=&\,
\int_t^T \frac{\langle\alpha,\kappa^\prime(s)\rangle}
{\langle\alpha,\kappa(s)\rangle} e^{-i\langle\alpha,\kappa(s)\rangle
\langle\alpha,gy\rangle}ds \notag\\
=&\,\int_{\langle\alpha, gy\rangle \langle\alpha,\kappa(t)\rangle}^{\langle\alpha,
gy\rangle\langle\alpha,\kappa(T)\rangle}
\frac{1}{u}e^{-iu}du. \notag
\end{align}
This integral exists, because the admissible curve $\kappa$ remains in the Weyl chamber $C$.
For abbreviation, put $\,\varphi_{\alpha,g}(t):=\
\langle\alpha,\kappa(t)\rangle |\langle\alpha,gy\rangle|,$
which is strictly positive for all $t$. Notice also that $\,\lim_{t\to\infty}
\varphi_{\alpha,g}(t) = +\infty$,  by admissibility of $\kappa$.
We thus obtain
\begin{equation}\label{(1.21)}
\lim_{T\to\infty} \int_t^T
\bigl(B_\alpha^\kappa(s)\bigr)_{g,\sigma_\alpha g}\, ds \,=
\,i\,\text{sign}(\langle\alpha,gy\rangle) \text{si}
\big(\varphi_{\alpha,g}(t)\big)
 -\text{Ci}\big(\varphi_{\alpha,g}(t)\big)
\end{equation}
where for $\tau >0$,  $\, \text{si}(\tau) =
- \int_\tau^\infty \frac{\sin u}{u}du\,$ and $\,\text{Ci}(\tau) = -\int_\tau^\infty \frac{\cos u}{u}du\,$ are the integral sine
and cosine respectively.
Thus in particular,
 condition (1) is satisfied. To verify condition (2), notice first
that
 the matrix entries of  $\, A^\kappa(t) \widetilde A^\kappa(t)$
 are linear combinations with constant coefficients,
of terms of the following kind:
\[ I_{\alpha,\beta,g}(t)\,:=\,
\frac{\langle\alpha,\kappa^\prime(t)\rangle}{\langle\alpha,\kappa(t)\rangle}\,
e^{-i\langle\alpha,\kappa(t)\rangle\langle\alpha,gy\rangle}
\bigl(\pm i\, \text{si}\big(\varphi_{\beta,\sigma_\alpha g}(t)\big)
  - \text{Ci}\big(\varphi_{\beta,\sigma_\alpha g}(t)\big)\bigr), \]
with $g\in G, \,
\alpha,\,\beta\in R_+.$
In order to estimate the integral sine and cosine terms,
 we use that for $\tau >0$,
\begin{equation} \label{(1.2a)}
|\text{si}(\tau)| \,\leq 2/\tau\,,\quad |\text{Ci}(\tau)| \,\leq 2/\tau\,.
\end{equation}
In fact, integration by parts yields
\[ \si(\tau) \,=\, - \frac{\cos \tau}{\tau} +\, \int_\tau^\infty
\frac{\cos u}{u^2} du,\]
which readily implies the first part of \eqref{(1.2a)}, and the second one
is seen in a similar way.
Moreover, as $\kappa$ is admissible, we have
$\langle\alpha,\kappa(t)\rangle >0$ for
all $t>0$, and there exists some constant
$\delta >0$ such that
\[\langle\beta,\kappa(t)\rangle\geq\, \delta |\kappa(t)|\,\geq\,
\delta\langle\alpha,\kappa(t)\rangle/\sqrt 2\,\quad\text{ for all }\,
\alpha,\beta\in R_+.\]
Together with \eqref{(1.2a)}, this yields the estimation
\[\big|I_{\alpha, \beta, g}(t)\big|\,\leq\,
C_1\frac{\langle\alpha,\kappa^\prime(t)\rangle}
{\langle\alpha,\kappa(t)\rangle\,\langle\beta,\kappa(t)\rangle}\, \leq\,
C_2\frac{\langle\alpha,\kappa^\prime(t)\rangle}
{\langle\alpha,\kappa(t)\rangle^2}\,
\]
with constants $C_1,\, C_2$ depending on $g, \alpha,\beta$ only.
As
\[ \int_{t_0}^\infty \frac{\langle\alpha,\kappa^\prime(t)\rangle}
{\langle\alpha,\kappa(t)\rangle^2}\,dt\, =\,
\int_{\langle\alpha,\kappa(t_0)\rangle}^\infty
\frac{1}{u^2}\,du\, <\infty,\]
it follows that $A^\kappa$ fulfills condition (2).
\end{proof}

\section{Proofs of the main theorems}

For the proof of Theorem \ref{T:main1} we consider $x\mapsto F(x,y)$ along arbitrary
admissible curves, and infer by an interpolation technique  that the limit in Theorem
\ref{T:Main} is independent of the special choice of the admissible curve $\kappa$ and also
of $y\in C$. We start with a supplementary notation.

\begin{definition}
A sequence $(x_n)_{n\in \b N} \subset C$ with $\lim_{n\to\infty} x_n \,=\, \infty$ is
called admissible, if there exists an interpolating admissible curve for it, i.e. an
admissible   $\,\kappa:(0,\infty)\to C$ such that $x_n = \kappa(t_n)$ for
suitable parameters $t_n$ with $\,\lim_{n\to\infty}t_n = \infty$.
\end{definition}

\begin{remark}
The following special situation will be of importance in the sequel:
Suppose that $(x_n)_{n\in \b N}$ is  contained in  $C_\delta$ for some
$\delta >0$
and satisfies $\lim_{n\to\infty}|x_n|=\infty$ as well as
$\,x_{n+1} - x_n \in C\,$ for all $n\in \b N$. Then $(x_n)_{n\in \b N}$ is
admissible. An admissible interpolating curve is obtained by
slightly smoothening the
piecewise linear connection of the successive points $x_n$.
\end{remark}

\begin{proof}[Proof of Theorem \ref{T:main1}]
In a first step, we show that there exists a non-zero  vector
  $v(y) = (v_g(y))_{g\in G}\in \b C^{|G|}$ such that for each admissible curve $\kappa$ in $C$,
\begin{equation}\label{(1.3)} \lim_{t\to\infty} F(\kappa(t),y)\,=\,v(y).
\end{equation}
For this, fix $y\in C$ and let $\kappa_1\,,\, \kappa_2$ be any two admissible curves, both
 contained in some $C_\delta$. With the above remark in mind, we can inductively construct an admissible sequence
 $(x_n)_{n\in \b N} \subset C_\delta$ with $x_{2n-1}\in \kappa_1$ and $x_{2n}\in \kappa_2$
for all $n\in \b N$. In fact, suppose that $x_1,\ldots, x_n$ are already
constructed, and consider the set $S_n=\{x\in C_\delta: x-x_n\notin C\}$, which is bounded. The curves $\kappa_i$ being admissible, we can therefore choose $x_{n+1}$ from the part of the relevant curve $\kappa_i$ which is contained in the complement of $S_n$, and we can do this in such a way that $\lim_{n\to\infty} |x_n|=\infty$.
 Now join the successive points $x_n$  by an interpolating
admissible curve $\kappa$. Then
according to Theorem \ref{T:Main}, all three limits
\[ \lim_{t\to\infty} F(\kappa_1(t),y),\quad \lim_{t\to\infty}
F(\kappa_2(t),y),
    \quad \lim_{t\to\infty} F(\kappa(t),y)\]
exist, are different from zero, and must in fact be equal by our choice of the
interpolating curve $\kappa$.
This proves \eqref{(1.3)}.
\noindent
Next, we focus on admissible rays. Observe that $\,F_g(tx,y) = F_{g^{-1}}(ty,x)$ for all $g\in G$ and
$x,y\in C$. Together with
\eqref{(1.3)}, this implies that
 $\, v_g(y) = v_{g^{-1}}(x)$, and therefore
also
 $\,v_g(x)\, = \,v_{g^{-1}}(x)\,=\, v_g(y)\, =: v_g.$  Put $v= (v_g)_{g\in G}$. Then
\begin{equation}\label{(4.10)} \lim_{t\to\infty} F(\kappa(t),y) = v
\end{equation}
for every admissible $\kappa$ and every $y\in C$. Now assume that the statement of Theorem 1 is
false. Then there exist  $\epsilon >0$ and
 a sequence $(x_n)_{n\in \b N} \subset C_\delta$
with $\,\lim_{n\to\infty} |x_n| = \infty$ and such that
\[ \max_{g\in G}|F_g(x_n,y) - v_g| >\epsilon \quad\text{for all }\,n\in \b N.\]
We may also assume without restriction that $(x_n)_{n\in \b N}$ is admissible - again because
for each $x\in C_\delta$ the set $\{z\in C_\delta: z-x\notin C\}$ is bounded. Hence relation
\eqref{(4.10)}
 entails $\,\lim_{n\to\infty}F_g(x_n,y) = v_g$, a contradiction.
 \end{proof}

The proof of Theorem \ref{T:main3} is based on the Phragm\'en-Lindel\"of
Theorems (see  Section 5.6. of \cite{T}) for
the right half plane $H=\{z\in \b C: \text{Re}\, z\geq 0\}$.

\begin{proof}[Proof of Theorem \ref{T:main3}]
We may assume that $\gamma>0$.
For fixed $x,y\in C$  define
\[ G(z) := z^\gamma\,\sqrt{w_k(x)w_k(y)}\, e^{-z\langle x,y\rangle}E_k(zx,y),\]
which is regular in $H\setminus\{0\}$ and continuous in $H$ with $G(0)=0$.
The integral representation \eqref{(2.77)} easily implies that
\[ |E_k(zx,y)|\,\leq \,\text{max}_{g\in G}\, e^{\text{Re}z \langle gx,y\rangle}
\quad\text{for all } z\in \b C, \]
c.f. Corollary 5.4 of \cite{R} or, alternatively, \cite{dJ1}. As $x$ and $y$ are contained
 in the same Weyl chamber, the inequality
$\,\langle gx,y\rangle \,\leq \langle x,y\rangle\,$ holds for all $g\in G$
(Theorem 3.1.2 of \cite{GB}). This shows that
\[ |G(z)|\leq\, |z|^\gamma\sqrt{w_k(x)w_k(y)}\,
 e^{\text{Re}z(\langle gx,y\rangle - \langle x,y\rangle)}\,
\leq \,|z|^\gamma\sqrt{w_k(x)w_k(y)}\,\]
as long as $\,\text{Re}\, z\geq 0$. Hence $G$ is of subexponential growth when restricted to $H.$
More precisely: for every  $\delta >0,$
\[ G(z) = O\bigl(e^{\delta |z|}\bigr) \quad\text{as } z\to\infty \,\text{ within }\,H.\]
Next consider $G$ along the boundary lines of $H,$
$\kappa_{\pm i}(t) = \pm it,\, t>0$.  According to Theorem \ref{T:main1},
$\,\lim_{t\to\infty} G(it) = i^\gamma v_e$. Moreover,  $\,G(-it) = \overline{G(it)}$
 for  $t>0$ (c.f. \eqref{(2.72)}); hence $\lim_{t\to\infty} G(-it)$
exists as well.
Employing the
Phragm\'en-Lindel\"of Theorems 5.62 and 5.64 of \cite{T}, we deduce that
$G$ is in fact  bounded in $H$
and that $\lim_{z\to\infty,\, z\in H} G(z) =  i^\gamma v_e\,.$
\end{proof}

We finally come to the proof of Theorem \ref{T:main2}. The key for
our approach
 is the following
simple observation: according to
formula \eqref{(2.77)}, one may write
\begin{equation}\label{(4.1)}
 E_k(x,-i\xi) =\, \int_{\b R^N}
e^{-i\langle \xi, y\rangle}d\mu_x^k(y)\,=\, \widehat{\mu_x^k}(\xi) \quad
(x,\xi\in \b R^N),
\end{equation}
where $\widehat \mu$ stands for the (classical)
Fourier-Stieltjes transform of  $\mu\in M^1(\b R^N)$,
\[ \widehat\mu(\xi) \,=\, \int_{\b R^N} e^{-i\langle \xi,y\rangle} d\mu(y).\]
Equation \eqref{(4.1)} suggests to employ Wiener's theorem,
which characterizes Fourier-Stieltjes transforms
of continuous measures on
locally compact abelian groups (here $(\b R^N, +)$),
 see for instance Lemma 8.3.7 of \cite{GMG}:

\begin{lemma} (Wiener) \label{L:Wiener}
For  $\mu\in M^1(\b R^N)$ the following
properties are equivalent: \parskip=-2pt
\begin{enumerate}\itemsep=-1pt
\item[\rm{(1)}] $\mu$ is continuous.
\item[\rm{(2)}] $\displaystyle
\lim_{n\to\infty} \frac{1}{n^N} \int_{\{\xi\in \b R^N:\,|\xi|\leq n\}} |\widehat \mu(\xi)|^2
d\xi\,=\, 0.$
\end{enumerate}
\end{lemma}

Apart from this, our argumentation relies on  the following growth estimate for
$E_k$, which is an easy
consequence of Theorem \ref{T:main1} and of some interest in its own:

\begin{proposition}
\label{P:Bound}
 Let $y\in C$. Then for each $\delta>0$
there exists a constant $M_\delta(y) >0$ such that
\[w_k(x)|E_k(ix,gy)|^2\,\leq M_\delta(y)
\quad\text{for all }\,x\in C_\delta\,,\, g\in G\]
\end{proposition}

\begin{remark}
It is important  at this point to notice that  the asymptotics
 of Theorem \ref{T:Main} implies boundedness of
$ x\mapsto w_k(x)|E_k(ix,gy)|^2 \,$ only within suitable subsets of $C$.
We do not know at present whether this function remains bounded when the range
of $x$ is all of $C$.
\end{remark}

\begin{proof}[Proof of Theorem \ref{T:main2}] For $n\in \b N$ put
$\,K_n:= \{\xi \in \b R^N: 1\leq |\xi|
\leq n\}.\,$
In view of the  properties of $E_k$ (\eqref{(2.71)} and  \eqref{(2.72)}),
 the assertion of the theorem
is equivalent to
\begin{equation}\label{(4.2)}
\lim_{n\to\infty} \frac{1}{n^N} \int_{K_n\cap C}
|E_k(ix,\xi)|^2 d\xi\,\,=\, 0
\quad\text{ for all }\, x\in \b R^N_{\,reg}.
\end{equation}
For fixed $x\in \b R^N_{\,reg}$ and $\delta >0$, define

\begin{align} I^\delta_1(n):=& \frac{1}{n^N} \int_{K_n\cap C_\delta} |E_k(ix,\xi)|^2 d\xi,\notag\\
I^\delta_2(n):= &\frac{1}{n^N} \int_{K_n\cap (C\setminus C_\delta)}
|E_k(ix,\xi)|^2 d\xi.\notag
\end{align}
Let further $\omega$ denote the
Lebesgue surface measure on $S^{N-1} = \{\xi \in \b R^N: |\xi|=1\}$.
By use of \eqref{(2.72)}, one obtains
\begin{align}
 I_2^\delta(n) \,&=\,\frac{1}{n^N} \int_1^n \int_{S^{N-1}
\cap(C\setminus C_\delta)} |E_k(ix,t\xi)|^2 d\omega(\xi)\, t^{N-1}dt
\notag\\
&\leq\, \frac{\omega\bigl(S^{N-1}\cap (C\setminus C_\delta)\bigr)}{n^N}
 \int_1^n t^{N-1} dt \,\, \leq \, \frac{1}{N} \cdot
\omega\bigl(S^{N-1}\cap (C\setminus C_\delta)\bigr), \notag
\end{align}
which tends to $0$ as $\delta\to 0$.
Thus for given $\epsilon >0$, we can find
 $\delta >0$
such that $\, I_2^\delta(n) \leq \epsilon \,$ for all $n\in \b N$.
With  this $\delta$ fixed,  the upper estimate on $E_k$ according to
Proposition \ref{P:Bound} yields
\[ I_1^\delta(n)\,\leq\,
   \frac{M_\delta(x)}{n^N} \int_{K_n\cap C_\delta} \frac{d\xi}{w_k(\xi)}
   \,.\]
The weight   $w_k$ being homogeneous of
degree $2\gamma$,
we further have
\[ \int_{K_n\cap C_\delta} \frac{d\xi}{w_k(\xi)}\,=\, A_\delta
\int_1^n t^{N-2\gamma-1} dt,\]
with
\[ A_\delta = \int_{S^{N-1}\cap C_\delta}
\frac{d\omega(\xi)}{w_k(\xi)}\, <\infty.\]
As $\gamma$ is strictly positive,  this implies that
$\, \lim_{n\to\infty} I_1^\delta(n) \,=\, 0$
and finishes the proof.
\end{proof}

\begin{remark} As already mentioned, we conjecture  that
for $x\in \b R^N_{reg}$ the measures $\mu_x^k$ are even absolutely
continuous with respect to Lebesgue measure, provided $R^\prime:=\{\alpha\in R\mid k(\alpha)> 0\}$ spans $\b R^N$.
We comment briefly on the hypotheses for this conjecture. First notice that some regularity condition on $x$ is necessary: In fact, for $x=0$
the representing measure is always given by the unit mass at the origin.
 As to the condition on $R^\prime$, let $V^\prime$ be the span of $R^\prime$, where $V^\prime=\{0\}$ by convention if $R^\prime=\emptyset$.
Suppose $V^\prime\not=\b R^n$. Let $V^{\prime\prime}:= (V^\prime)^\bot$ be the orthoplement with corresponding decomposition $\b R^N = V^\prime \oplus
V^{\prime\prime}$. For $x\in \b R^N$, write $x= x^\prime + x^{\prime\prime}$ with
$x^\prime\in V^\prime\,, \, x^{\prime\prime}\in V^{\prime\prime}.$ Then it is easily seen from the characterization \eqref{(2.66)}
of $E_k$  that
\[ E_k(x,y) \,=\,  E_k(x^{\prime},y^{\prime})\cdot
e^{\langle x^{\prime\prime},y^{\prime\prime}
\rangle} \quad \text{for all }\, x,y\in \b R^N\]
and accordingly,
\[ \mu_x^k\, =\, \mu_{x^\prime}^k \otimes
\delta_{x^{\prime\prime}},\]
where on the right side, $k$ is understood as a multiplicity function on $R^\prime$.
Thus for all $x\in \b R^N$, $\mu_x^k$ is supported in a translate of $V^\prime$ and is therefore not absolutely continuous.

\end{remark}

\begin{acknowledge}
For the first author it is a pleasure to thank the Isaac Newton
Institute in Cambridge, U.K. for their hospitality during the preparation
of  this article.
\end{acknowledge}

\end{document}